\documentclass [11pt,twoside]   {article}

  \newtheorem{theorem}{Theorem}
  \newtheorem{lemma}{Lemma}
  \newtheorem{proposition}{Proposition}

  \newtheorem{remark}{Remark}

  \newcommand {\pf}  {\mbox{\sc Proof. \,\,}}

  \newcommand {\qed} {\null \hfill \rule{2mm}{2mm}}

\begin {document}

\title{{\Large{\bf Bounds for the size of sets with the property $D(n)$ }}}

\author
{{\sc Andrej Dujella} \vspace{1ex}\\ University of Zagreb,
Croatia}

\date{}
\maketitle

\begin{abstract}
\noindent Let $n$ be a nonzero integer and $a_1<a_2<\cdots<a_m$
positive integers such that $a_ia_j+n$ is a perfect square for all
$1\leq i<j\leq m$. It is known that $m\leq 5$ for $n=1$. In this
paper we prove that $m\leq 31$ for $|n|\leq 400$ and $m<15.476\,
\log{|n|}$ for $|n|>400$.
\end{abstract}

\footnotetext{ {\it 2000 Mathematics Subject Classification.}
11D45, 11D09, 11N36.

{\it Key words and phrases.} Diophantine $m$-tuples, property
$D(n)$, large sieve.}

\section{Introduction}
Let $n$ be a nonzero integer. A set of $m$ positive integers
$\{a_1,a_2,\ldots,a_m\}$ is called \emph{a $D(n)$-$m$-tuple} (or
\emph{a Diophantine $m$-tuple with the property $D(n)$}) if
$a_ia_j+n$ is a perfect square for all $1\leq i<j\leq m$.

Diophantus himself found the $D(256)$-quadruple
$\{1,\,33,\,68,\,105\}$, while the first $D(1)$-quadruple, the set
$\{1,\,3,\,8,\,120\}$, was found by Fermat (see \cite{Dic,Dio}).
In 1969, Baker and Davenport \cite{B-D} proved that this Fermat's
set cannot be extended to a $D(1)$-quintuple, and in 1998, Dujella
and Peth\H{o} \cite{D-P} proved that even the Diophantine pair
$\{1,3\}$ cannot be extended to a D(1)-quintuple. A famous
conjecture is that there does not exist a $D(1)$-quintuple. We
proved recently that there does not exist a $D(1)$-sextuple and
that there are only finitely many, effectively computable,
$D(1)$-quintuples (see \cite{D-jnt, D-fin}).

The question is what can be said about the size of sets with the
property $D(n)$ for $n\neq 1$. Let us mention that Gibbs
\cite{Gibbs1} found several examples of Diophantine sextuples,
e.g. $\{99,\,315,\,9920,\,32768,\,44460,\,19534284\}$ is a
$D(2985984)$-sextuple.

Define
\[ M_n=\sup \{ |S| \,:\, \mbox{$S$ has the property $D(n)$} \}. \]
Considering congruences modulo $4$, it is easy to prove that
$M_{n}=3$ if $n\equiv 2 \!\!\pmod{4}$ (see \cite{Bro,G-S,M-R}). On
the other hand, if $n\not\equiv 2\!\!\pmod{4}$ and
$n\not\in\{-4,\,-3,\,-1,\,3,\,5,\,8,\,12,\,20\}$, then $M_n\geq 4$
(see \cite{D-acta1}).

In \cite{D-size}, we proved that $M_n \leq 32$ for $|n|\leq 400$
and
\[ M_n < 267.81 \,\log{|n|} \,(\log\log{|n|})^2 \quad
\mbox{for $|n| > 400$}.\] The purpose of the present paper is to
improve this bound for $M_n$, specially in the case $|n| > 400$.
We will remove the factor $(\log\log{|n|})^2$, and also the
constants will be considerably smaller.

The above mentioned bounds for $M_n$ were obtained in
\cite{D-size} by considering separately three types (large, small
and very small) of elements in a $D(n)$-$m$-tuple. More precisely,
let
\begin{eqnarray*}
A_n \!\!&=&\!\! \sup \{ |S\cap [|n|^3, +\infty \rangle |\,:\,
\mbox{$S$ has the property $D(n)$} \}, \\
B_n \!\!&=&\!\! \sup \{ |S\cap \langle n^2, |n|^3 \rangle |\,:\,
\mbox{$S$ has the property $D(n)$} \}, \\
C_n \!\!&=&\!\! \sup \{ |S\cap [1,n^2] |\,:\,
\mbox{$S$ has the property $D(n)$} \}.
\end{eqnarray*}

In \cite{D-size}, it was proved that $A_n\leq 21$ and $B_n < 0.65
\,\log{|n|} +2.24$ for all nonzero integers $n$, while $C_n <
265.55 \,\log{|n|}\,(\log\log{|n|})^2 +9.01 \,\log\log{|n|}$ for
$|n| > 400$ and $C_n \leq 5$ for $|n|\leq 400$. The combination of
these estimates gave the bound for $M_n$.

In the estimate for $A_n$, a theorem of Bennett \cite{Ben} on
simultaneous approximations of algebraic numbers was used in
combination with a gap principle, while a variant of the gap
principle gave the estimate for $B_n$. The bound for $C_n$ (number
of "very small" elements) was obtained using the Gallagher's large
sieve method \cite{Gal} and an estimate for sums of characters.

In the present paper, we will significantly improve the bound for
$C_n$ using a result of Vinogradov on double sums of Legendre's
symbols. Let us mention that Vinogradov's result, in a slightly
weaker form, was used recently, in similar context, by Gyarmati
\cite{Gya} and S\'ark\"ozy \& Stewart \cite{Sar-Ste}. We will
prove the following estimates for $C_n$.

\begin{proposition} \label{pr:1}
If $|n| > 400$, then \,$C_n < 11.006 \,\log{|n|}$. If $|n|\geq
{10}^{100}$, then \,$C_n < 8.37 \,\log{|n|}$.
\end{proposition}

More detailed analysis of the gap principle used in \cite{D-size}
will lead us to the slightly improved bounds for $B_n$.

\begin{proposition} \label{pr:2}
For all nonzero integers $n$ it holds \,$B_n < 0.6114\, \log{|n|}
+ 2.158$. If $|n| > 400$, then \,$B_n < 0.6071\, \log{|n|} +
2.152$.
\end{proposition}

By combining Propositions \ref{pr:1} and \ref{pr:2} with the above
mentioned estimate for $A_n$, we obtain immediately the following
estimates for $M_n$.

\begin{theorem} \label{tm:1}
If $|n| \leq 400$, then $M_n \leq 31$. If $|n| > 400$, then \,$M_n
< 15.476, \log{|n|}$. If $|n|\geq {10}^{100}$, then \,$M_n <
9.078\, \log{|n|}$.
\end{theorem}

\section{Three lemmas}

\begin{lemma}[Vinogradov] \label{l:V}
Let $p$ be an odd prime and $\gcd(n,p)=1$. If $A,B \subseteq
\{0,1, \ldots, p-1\}$ and
\[ T = \sum_{x\in A} \sum_{y\in B} \Big( \frac{xy+n}{p} \Big), \]
then $|T| < \sqrt{p|A|\cdot |B|}$.
\end{lemma}

\pf See \cite[Problem V.8.c)]{Vin}. \qed

\begin{lemma}[Gallagher] \label{l:G}
If all but $g(p)$ residue classes {\rm mod} $p$ are removed for
each prime $p$ in a finite set ${\cal S}$, then the number of
integers which remain in any interval of length $N$ is at most
\begin{equation} \label{gal}
 \Big( \sum_{p\in{\cal S}} \log{p} -\log{N}\Big)
\Big/ \Big( \sum_{p\in{\cal S}} \frac{\log{p}}{g(p)} -
\log{N}\Big)
\end{equation}
provided the denominator is positive.
\end{lemma}

\pf See \cite{Gal}. \qed

\begin{lemma} \label{l:3}
If $\{a,b,c\}$ is a Diophantine triple with the property $D(n)$
and $ab+n=r^2$, $ac+n=s^2$, $bc+n=t^2$, then there exist integers
$e$, $x$, $y$, $z$ such that
\[ ae+n^2=x^2, \quad be+n^2=y^2, \quad ce+n^2=z^2 \]
and
\[ c=a+b+\frac{e}{n} + \frac{2}{n^2}(abe+rxy). \]
\end{lemma}

\pf See \cite[Lemma 3]{D-size}. \qed

\section{Proof of Proposition 1}

Let $N\geq n^2$ be a positive integer. Since $|n|>400$, we have
$N>1.6\cdot 10^5$. Let $D =\{a_1,a_2,\ldots,a_m\} \subseteq \{1,2,
\ldots, N\}$ be a Diophantine $m$-tuple with the property $D(n)$.
We would like to find an upper bound for $m$ in term of $N$. We
will use the Gallagher's sieve (Lemma \ref{l:G}). Let
\[ {\cal S}=\{ p \,:\, \mbox{$p$ is prime, $\gcd(n,p)=1$ and $p\leq Q$} \}, \]
where $Q$ is sufficiently large. For a prime $p\in {\cal S}$, let
$C$ denotes the set of integers $b$ such that $b \in
\{0,1,2,\dots, p-1\}$ and there is at least one $a\in D$ such that
$b\equiv a \!\!\pmod{p}$. Then $\Big( \frac{xy+n}{p} \Big) \in
\{0,1\}$ for all distinct $x,y\in C$. Here $\Big( \frac{.}{p}
\Big)$ denotes the Legendre symbol. If $0\in C$, then $\Big(
\frac{n}{p} \Big) =1$. For given $x\in C\setminus \{0\}$, we have
$\Big( \frac{xy_0+n}{p} \Big)$ for at most one $y_0\in C$. If
$y\neq x,y_0$, then $\Big( \frac{xy+n}{p} \Big)=1$. Therefore,
\begin{eqnarray*}
T &=& \sum_{x,y\in C} \Big( \frac{xy+n}{p} \Big) = \sum_{x\in C}
\Big( \sum_{y\in C} \Big( \frac{xy+n}{p} \Big) \Big) \\ &\geq &
\sum_{x\in C} (|C|-3) \geq |C|(|C|-3).
\end{eqnarray*}
On the other hand, Lemma \ref{l:V} implies
\[ T < |C| \cdot \sqrt{p}. \]
Thus, $|C| < \sqrt{p} +3$ and we may apply Lemma \ref{l:G} with
\[ g(p) = \min \{ \lfloor \sqrt{p} \rfloor +3, p\}. \]

Let us denote the numerator and denominator from (\ref{gal}) by
$E$ and $F$, respectively. By \cite[Theorem 9]{R-S}, we have
\[ E=\sum_{p\in {\cal S}} \log{p} - \log N<\theta(Q) <
1.01624\,Q.\] The function $f(x)=\frac{\log x}{\min
\{\sqrt{x}+3,x\}}$ is strictly decreasing for $x>25$. Also, if
$Q\geq 118$, then $f(p)\geq f(Q)$ for all $p\leq Q$.

For $p\in \mathcal{S}$ it holds $\gcd(n,p)=1$. This condition
comes from the assumptions of Lemma \ref{l:V}. However, we will
show later that $n$ can be divisible only by small proportion of
the primes $\leq Q$. Assume that $n$ is divisible by at most 5\%
of primes $\leq Q$.
Then, for $Q\geq 118$, we have
\begin{eqnarray} \label{z2}
F &\geq& \sum_{p\in{\cal S}} f(p) - \log N \geq \frac{\log
Q}{\sqrt{Q} +3} \cdot |{\cal S}| - \log{N} \nonumber \\ &\geq&
\frac{\log Q}{\sqrt{Q} +3} \cdot \frac{19}{20} \pi(Q) - \log{N} >
\frac{0.95\, Q}{\sqrt{Q}+3} - \log N.
\end{eqnarray}
Since $F$ has to be positive in the applications of Lemma
\ref{l:G}, we will choose $Q$ of the form
\begin{eqnarray*}
Q= c_1 \cdot \log^2{N}.
\end{eqnarray*}

We have to check whether our assumption on the proportion of
primes which divide $n$ is correct. Suppose that $n$ is divisible
by at least 5\% of the primes $\leq Q$. Then $|n|\geq p_1p_2\cdots
p_{\lceil \pi(Q)/20 \rceil}$, where $p_i$ denotes the $i$-th
prime. By \cite[3.5 and 3.12]{R-S}, we have $p_{\lceil \pi(Q)/20
\rceil} > R$, where
\[ R= \frac{1}{20}\frac{Q}{\log{Q}}
\log\Big(\frac{1}{20}\frac{Q}{\log{Q}} \Big).\] Assume that
$c_1\geq 6$. Then $Q>860$ and $R>11.77$. From \cite[3.16]{R-S}, it
follows that
\begin{equation} \label{z1}
 \log{|n|}>\sum_{p\leq R} \log{p} >R\Big(1-\frac{1.136}{\log{R}}
\Big).
\end{equation}
Furthermore, $\frac{1}{20}\frac{Q}{\log{Q}}
>Q^{0.273}$ and $R>0.0136\,Q$. Hence, (\ref{z1}) implies
$\log{R}>7.793$ and therefore
\[ \log N \geq 2\log |n| > 0.01466\, Q \geq  0.08796\, \log^2 N, \]
contradiction the assumption that $N>1.6\cdot 10^5$.

Therefore, we have that $n$ is divisible by at most 5\% of the
primes $\leq Q$, and hence we have justifies the estimate
(\ref{z2}).

Under the assumption that $c_1\geq 6$, the inequality (\ref{z2})
implies
\[ F> 0.861\, \sqrt{Q} - \log N = (0.861\,\sqrt{c_1} -1) \log N \]
and
\[ \frac{E}{F} < \frac{1.017\, c_1}{0.861\,\sqrt{c_1}-1} \cdot \log
N.\] For $c_1=6$ we obtain
\begin{equation} \label{400}
\frac{E}{F} < 5.503\, \log N.
\end{equation}

Assume now that $N\geq {10}^{200}$ and $c_1\geq 4$. Then $Q>
848303$ and we can prove in the same manner as above that $n$ is
divisible by at most 1\% of the primes $\leq Q$. This fact implies
\[ \frac{E}{F} < \frac{1.017 c_1}{0.986\sqrt{c_1}-1} \cdot \log
N.\] For $c_1=4.11$ we obtain
\begin{equation} \label{100}
\frac{E}{F} < 4.185 \log N.
\end{equation}

Setting $N=n^2$ in (\ref{400}) and (\ref{100}), we obtain the
statements of Proposition \ref{pr:1}. \qed

\section{Proof of Proposition 2}

We may assume that $|n| > 1$. Let $\{a,b,c,d\}$ be a
$D(n)$-quadruple such that $n^2 < a<b<c<d$. We apply Lemma
\ref{l:3} on the triple $\{b,c,d\}$. Since $b>n^2$ and $be+n^2
\geq 0$, we have that $e\geq 0$. If $e=0$, then
$d=b+c+2\sqrt{bc+n} < 2c+2\sqrt{c(c-1)+n} < 4c$, contradicting the
fact that $d>4.89\, c$ (see \cite[Lemma 5]{D-size}).

Hence $e\geq 1$ and
\begin{equation} \label{d}
d > b+c+\frac{2bc}{n^2}+
\frac{2t\sqrt{bc}}{n^2}.
\end{equation}
Lemma \ref{l:3} also implies
\begin{equation} \label{c}
c \geq a+b+2r.
\end{equation}
From $r^2\geq ab-\sqrt[4]{ab}$ and $ab\geq 30$ it follows that $r
> 0.96\, a$, and (\ref{c}) implies $c>3.92\, a$. Similarly, $bc\geq 42$ implies
$t > 0.969\, \sqrt{bc}$ and, by (\ref{d}), $d> b+c+ 3.938\,
\frac{bc}{n^2}
> 4.938\, c+b$.

Assume now that $\{a_1,a_2,\ldots,a_m\}$ is a $D(n)$-$m$-tuple and
$n^2 < a_1 < a_2 < \cdots <a_m <|n|^3$. We have
\begin{eqnarray*}
a_3 > 3.92\, a_1, \quad a_{i} > 4.938\, a_{i-1} + a_{i-2}, \quad
\mbox{for $i=4,5,\ldots,m$}.
\end{eqnarray*}
Therefore, $a_m > {\alpha}_m a_1$, where the sequence
$({\alpha}_k)$ is defined by
\begin{eqnarray} \label{alpha}
{\alpha}_k = 4.938 {\alpha}_{k-1}+{\alpha}_{k-2}, \quad {\alpha}_2
= 1, \,\, {\alpha}_3 = 3.92.
\end{eqnarray}
Solving the recurrence (\ref{alpha}), we obtain ${\alpha}_k
\approx \beta {\gamma}^{k-3}$, with $\beta \approx 3.964355$,
$\gamma \approx 5.132825$. More precisely,
\[ |{\alpha}_k - \beta {\gamma}^{k-3}| < \frac{1}{\beta
{\gamma}^{k-3}} \,.\] From $|n|^3 - 1 \geq  a_m > {\alpha}_m a_1
\geq {\alpha}_m (n^2+1)$, it follows ${\alpha}_m \leq
|n|-\frac{1}{|n|}$ and $\beta {\gamma}^{m-3} < |n|$. Hence,
\begin{equation} \label{m}
m < \frac{1}{\log{\gamma}} \log{|n|} + 3 -
\frac{\log{\beta}}{\log{\gamma}}.
\end{equation}

For the above values of $\beta$ and $\gamma$ we obtain
\[ m < 0.6114\, \log{|n|} + 2.158. \]

Assume now that $|n| > 400$. Then $bc>ab>{400}^4$, which implies
$c>3.999999\, a$ and $d>4.999999\, c + b$. Therefore, in this case
the relation (\ref{m}) holds with $\beta \approx 4.042648$,
$\gamma \approx 5.192581$, and we obtain
\[ m < 0.6071\, \log{|n|} + 2.152. \]
\qed

\begin{remark}
{\rm The constants in Theorem \ref{tm:1} can be improved, for
large $|n|$, by using formula (2.26) from \cite{R-S} in the
estimate for the sum $\sum_{p \in {\cal S}} f(p)$. In that way, it
can be proved that for every $\varepsilon > 0$, $F >
(2-\varepsilon) \sqrt{Q} - \log N$ holds for sufficiently large
$Q$.

Also, in the proof of Proposition \ref{pr:2}, for sufficiently
large $|n|$ we have $c>(4-\varepsilon)a$ and
$d>(5-\varepsilon)c+b$, which leads to $\displaystyle{ B_n <
\bigg(\frac{1}{\log(\frac{5+\sqrt{29}}{2})}+\varepsilon \bigg)
\,\log |n|}$.

These results imply that for every $\varepsilon
> 0$ there exists $n(\varepsilon)$ such that for
$|n|>n(\varepsilon)$ it holds
\[ M_n < \bigg(2+\frac{1}{\log(\frac{5+\sqrt{29}}{2})}+\varepsilon \bigg)
\,\log |n|. \]
}
\end{remark}

\medskip {\bf Acknowledgements.} The author is grateful to
the referees for valuable comments, and in particular for pointing
out a gap in the proof of Proposition \ref{pr:1} in the first
version of the manuscript.

\medskip

{\small \noindent Department of Mathematics \\ University of
Zagreb
\\ Bijeni\v cka cesta 30, 10000 Zagreb \\ Croatia \\
{\em E-mail address}: {\tt duje@math.hr}}

\end{document}